\documentclass[a4paper,12pt,intlimits,oneside]{amsart}

\usepackage{enumerate}

\usepackage{latexsym,amssymb}

\newcommand{\comment}[1]{}

\newcommand{\epf}{ $\Box$\medskip}

\numberwithin{equation}{section}

\def\lsim{\raisebox{-1ex}{$~\stackrel{\textstyle <}{\sim}~$}}

\newcounter{rea}
\newcounter{rej}
\newcounter{res}
\newtheorem{thm}{Theorem}[section]
\newtheorem{prop}[thm]{Proposition}

\newtheorem{lem}[thm]{Lemma}

\begin{document}

\title[norm inequalities in generalized Morrey spaces]{Norm inequalities in generalized Morrey spaces}% a class of function spaces including Weighted Morrey spaces}
%\author[A. Bonami]{Aline Bonami}
%\address{MAPMO-UMR 6628,
%D\'epartement de Math\'ematiques, Universit\'e d'Orleans,
%45067 Orl\'eans Cedex 2, France}
%\email{{\tt Aline.Bonami@univ-orleans.fr}}
\author[J. Feuto]{Justin Feuto}
\address{Laboratoire de Math\'ematiques Fondamentales, UFR Math\'ematiques et Informatique, Universit\'e F\'elix Houphouet-Boigny (Cocody), 22 B.P 1194 Abidjan 22. C\^ote d'Ivoire}
\email{{\tt justfeuto@yahoo.fr}}

%43A15 ; 42B20 ; 42B25 ; 42B35 47B47}
\date{}

\begin{abstract}
 We prove that Calder\'on-Zygmund operators, Marcinkiewicz operators, maximal operators associated to Bochner-Riesz operators, operators with rough
kernel as well as commutators associated to these operators which are known to be bounded on weighted Morrey spaces under appropriate conditions, are bounded on a wide family of function spaces.% under the same hypothesis.
\end{abstract}

%\noindent {\bf AMS Subject Classification:} 

%\vspace{.08in} \noindent \textbf{Keywords}:
\maketitle 

\section{Introduction}
In the last decade, many works in classical harmonic analysis have been devoted to norm inequalities involving classical and non-classical operators in the setting of weighted Morrey spaces. The results obtained are most of the time extensions of well known  analogues in the weighted Lebesgue spaces. 

We equip the $n$-dimensional Euclidean space $\mathbb R^{n}$ with the Euclidean norm $\left|\cdot\right|$ and the Lebesgue measure $dx$. A weight is any positive measurable function $w$ which is locally integrable on $\mathbb R^{n}$. Let $w$ be a weight, $1\leq q<\infty$ and $0<\kappa<1$. As Komori and Shirai in \cite{KS}, we define the weighted Morrey space by
$$L^{q,\kappa}(w)=\left\{f\in L^{q}_{loc}(w):\left\|f\right\|_{L^{q,\kappa}(w)}<\infty\right\}$$
% consists of measurable functions $f$ such that $\left\|f\right\|_{L^{q,\kappa}(w)}<\infty$,
 where
\begin{equation}
\left\|f\right\|_{L^{q,\kappa}(w)}:=\sup_{B}\left(\frac{1}{w(B)^{\kappa}}\int_{B}\left|f(x)\right|^{q}w(x)dx\right)^{\frac{1}{q}}.\label{wmorrey}
\end{equation}
The supremum is taken over all balls $B$ in $\mathbb R^{n}$, and $w(B)=\int_{B}w(x)dx$. When $w\equiv 1$, we use the notation $\left|B\right|$ for the Lebesgue measure of $B$. These spaces can be viewed as extensions of weighted Lebesgue spaces $L^{q}(w)$,i.e., spaces that consist of measurable functions $f$ satisfying % We recall that $f\in L^{q}(w)$ if %is the spaces of Lebesgue measurable function $f$ such that 
\begin{equation*}
\left\|f\right\|_{L^{q}(w)}:=\left(\int_{\mathbb R^{n}}\left|f(x)\right|^{q}w(x)dx\right)^{\frac{1}{q}}<\infty.
\end{equation*}

It has been proved by many authors (see \cite{KS},\cite{Wa2},\cite{Wa}) that most of the operators which are bounded on a weighted Lebesgue space are also bounded in an appropriate weighted Morrey space. We are going to prove that these results are valid on a larger family  of functional spaces including weighted Morrey spaces. 

Let $w$ be a weight and $1\leq q\leq\alpha\leq p\leq\infty$. We define the space $(L^{q}(w),L^{p})^{\alpha}:=(L^{q}(w),L^{p})^{\alpha}(\mathbb R^{n})$ as the space of all measurable functions $f$ satisfying $\left\|f\right\|_{(L^{q}(w),L^{p})^{\alpha}}<\infty$, where 
\begin{equation*}
\left\|f\right\|_{(L^{q}(w),L^{p})^{\alpha}}:=\sup_{r>0}\ _{r}\left\|f\right\|_{(L^{q}(w),L^{p})^{\alpha}}, 
\end{equation*}
with
\begin{equation}
\ _{r}\left\|f\right\|_{(L^{q}(w),L^{p})^{\alpha}}:=\left[\int_{\mathbb R^{n}}\left(w(B(y,r))^{\frac{1}{\alpha}-\frac{1}{q}-\frac{1}{p}}\left\|f\chi_{B(y,r)}\right\|_{L^{q}(w)}\right)^{p}dy\right]^{\frac{1}{p}}\label{normqp}
\end{equation}
for any $r>0$, with the usual modification when $p=\infty$.
 In the case $w\equiv 1$, we recover the space $(L^{q},L^{p})^{\alpha}$ defined in \cite{Fo0} by Fofana (see also \cite{FFK1}, \cite{FFK2}). Condition $q\leq\alpha\leq p$ ensures that the space is non trivial.  
For $q<\alpha$ and $p=\infty$, the space $(L^{q}(w),L^{\infty})^{\alpha}(\mathbb R^{n})$ is the weighted Morrey space $L^{q,\kappa}(w)$, with $\kappa=\frac{1}{q}-\frac{1}{\alpha}$. 

It is immediate that the spaces $(L^{q}(w),L^{p})^{\alpha}$ equipped with $\left\|\cdot\right\|_{(L^{q}(w),L^{p})^{\alpha}}$ are Banach spaces for $1\leq q\leq\alpha\leq p\leq\infty$. Let $1\leq q_{1}\leq q_{2}\leq\alpha\leq p\leq\infty$. For any weight $w$, we have the following inclusion : 
$$\left\|f\right\|_{(L^{q_{1}}(w),L^{p})^{\alpha}}\leq\left\|f\right\|_{(L^{q_{2}}(w),L^{p})^{\alpha}}$$
which comes from H\"older inequality. In the particular case where $w=1$, this family of spaces is an increasing family in $p$.i.e., $(L^{q},L^{p_{1}})^{\alpha}\subset (L^{q},L^{p_{2}})^{\alpha}$ whenever $1\leq q\leq \alpha\leq p_{1}\leq p_{2}$. To our knowledge, the similar inclusions for spaces with weight are still open problems. %But, up to now, we have no such result for the weighted spaces.  
 
 We prove in this paper that Calder\'on-Zygmund operators, Marcinkiewicz operators, maximal operators associated to Bochner-Riesz operators, operators with rough
kernel and associated commutators which are known to be bounded on weighted Morrey spaces under appropriate conditions, are bounded on this weighted spaces for appropriate weight. 
In fact, we prove that operators which are bounded on weighted Lebesgue spaces and satisfy some local pointwise control, are also bounded in our context. Operators fulfilling this conditions, include Littlewood Paley operators with rough kernels, whose control in this spaces was given by Wei and Tao in \cite{WT}. We point out that the paper of Wei and Tao was published while ours was already in the hands of the referee and available on arxiv.

 An important fact here is that the proof is simple and is the same for all kinds of operators that have been considered.

This paper is organized as follows. 

In the next section, we recall the definitions of the operators we are going to deal with, and state the main results. Section 3 is devoted to proofs. 

 Throughout the paper, the letter $C$ is used for non-negative constants independent of the relevant variables that may change from one occurrence to another. Constants with subscript, such as $C_{0}$, do not change in different occurrences. We propose the following abbreviation $\mathrm{\bf A}\lsim \mathrm{\bf B}$ for the inequalities $\mathrm{\bf A}\leq C\mathrm{\bf B}$, where $C$ is a positive constant independent of the main parameters. If we have $\mathrm{\bf A}\lsim \mathrm{\bf B}$ and  $\mathrm{\bf B}\lsim \mathrm{\bf A}$ then we put $\mathrm{\bf A}\cong \mathrm{\bf B}.$

 For $\lambda>0$ and a ball $B\subset\mathbb R^{n}$, we write $\lambda B$ for the ball with same center as $B$ and radius $\lambda$ times radius of $B$. We denote by  $E^{c}$ the complement of $E$. 
 
 \medskip

{\bf Acknowledgement.} The author is very grateful to Aline Bonami and the referee, for their useful comments and remarks.%  thank Aline Bonami for her support carefully reading and revision of the manuscript.
 %for the subset of $\mathbb R^{n}$ consisted of elements which are not in $E\subset\mathbb R^{n}$.% We denote by $\mathbb N^{\ast}$ the set of all positive integers.
 
\section{Definitions and statement of the main results}\label{section2}
 A weight $w$  belongs to $\mathcal A_{q}$ for $1\leq q<\infty$ if there exists a constant $C>0$ such that for all balls $B\subset\mathbb R^{n}$ we have
% for $q>1$%We recalled that $w\in\mathcal A_{p}$  if 
\begin{equation}
\left\{\begin{array}{lll}
\left(\frac{1}{\left|B\right|}\int_{B}w(x)dx\right)\left(\frac{1}{\left|B\right|}\int_{B}w^{\frac{-q'}{q}}(x)dx\right)^{\frac{q}{q'}}\leq C&\text{ if }&q>1,\\
\frac{1}{\left|B\right|}\int_{B}w(z)dz\leq C\mathrm{ess }\inf_{x\in B}w(x)&\text{ if }&q=1.
\end{array}\right.\label{apcondition}
\end{equation}
where $\frac{1}{q}+\frac{1}{q'}=1$.
%and for $q=1$
%\begin{equation}
%.
%\end{equation}
We put $\mathcal A_{\infty}=\cup_{q\geq 1}\mathcal A_{q}$. 

 It is known (see \cite{Gr} Proposition 9.1.5 p. 679 and Theorem 9.2.2 p. 685) that for $w\in\mathcal A_{p}$ with $1\leq p<\infty$:
 \begin{itemize}
 \item  the measure $w(x)dx$ is a doubling measure: precisely for all $\lambda>1$ and all balls $B$ we have 
 \begin{equation}
 w(\lambda B)\lsim \lambda^{np}w(B),\label{doubling}
 \end{equation}
 \item there exists $s>1$ such that for any ball $B\subset\mathbb R^{n}$, we have 
\begin{equation}
\left(\frac{1}{\left|B\right|}\int_{B}w(x)^{s}dx\right)^{\frac{1}{s}}\lsim\frac{1}{\left|B\right|}\int_{B}w(x)dx.\label{reverseholder}
\end{equation}
\end{itemize}
H\"older inequality and (\ref{reverseholder}) lead to  
\begin{equation}
\frac{w(E)}{w(B)}\lsim \left(\frac{\left|E\right|}{\left|B\right|}\right)^{\frac{s-1}{s}},\label{contwl}
\end{equation}
for any measurable subset $E$ of a ball $B$. 
%We say that $w$ satisfies a Gehring condition (or a reverse Hölder inequality) if there exists a $p>1$ %It is classic that for $1\leq q<\infty$, any weight $w$ belonging to $\mathcal A_{q}$ belongs to some reverse H\"older class $RH_{\tau}$ for $\tau>1$, that is there exists $C>0$ such that
%\begin{equation*}
%\left(\frac{1}{\left|B\right|}\int_{B}w(x)^{\tau}dz\right)^{\frac{1}{\tau}}\leq \frac{C}{\left|B\right|}\int_{B}w(x)dx.
%\end{equation*}
%Also, it is well known that if $w\in RH_{\tau}$ then there exists $C>0$ such that 
%\begin{equation}
%\frac{w(E)}{w(B)}\leq C\left(\frac{\left|E\right|}{\left|B\right|}\right)^{\frac{\tau-1}{\tau}},\label{rholder}
%\end{equation} 
%for any measurable subset $E$ of a ball $B$.

Let $T$ be a Calder\'on-Zygmund operator given by % i.e., there exists a function $K$ on $\mathbb R^{n}$ which satisfies the following conditions
\begin{equation*}
Tf(x)=\mathrm{p.v}\int_{\mathbb R^{n}}K(x-y)f(y)dy,
\end{equation*}
where $K$ is of class $\mathcal C^{1}(\mathbb R^{n}\setminus\left\{0\right\})$ with
\begin{equation*}
\left|K(x)\right|\leq\frac{C_{K}}{\left|x\right|^{n}}\text{ and }\left|\nabla K(x)\right|\leq\frac{C_{K}}{\left|x\right|^{n+1}}\text{ for }\ x\neq 0.
\end{equation*}
%As we can see ,if $1<q<\infty$ and $w\in\mathcal A_{q}$, then 
It is a classical result that the operator $T$ is bounded on $L^{q}(w)$ whenever $1<q<\infty$ and $w\in\mathcal A_{q}$, whereas for $q=1$ and $w\in\mathcal A_{1}$ we have the following weak type inequality % then for all $\lambda>0$, we have 
\begin{equation}
\left\|Tf\right\|_{L^{1,\infty}(w)}:=\sup_{\lambda>0}\lambda w(\left\{x\in\mathbb R^{n}:\left|Tf(x)\right|>\lambda\right\})\lsim\left\|f\right\|_{L^{1}(w)}.\label{wwleb}
\end{equation}
These results can be found in \cite{Du} or \cite{Gr}. Komori and Shirai extended them to weighted Morrey spaces in \cite{KS}. %, Theorem 3.3 ). 
%We can replace in Relation (\ref{normqp}) the weighted Lebesgue norm $\left\|\cdot\right\|_{q_{w}}$ by its weak version. That is,
We prove in this paper that the results remain valid in our weighted spaces. For this purpose, we put for $r>0$
\begin{equation*}
\ _{r}\left\|f\right\|_{(L^{1,\infty}(w),L^{p})^{\alpha}}:=\left[\int_{\mathbb R^{n}}\left(w(B(y,r))^{\frac{1}{\alpha}-1-\frac{1}{p}}\left\|f\chi_{B(y,r)}\right\|_{L^{1,\infty}(w)}\right)^{p}dy\right]^{\frac{1}{p}},
\end{equation*}
%where 
%\begin{equation}
%\left\|f\right\|^{\ast}_{q_{w},\infty}:=\sup_{\lambda>0}\lambda w(\left\{x\in\mathbb R^{n}/\left|f(x)\right|>\lambda\right\})^{\frac{1}{q}},
%\end{equation}
and
\begin{equation}
\left\|f\right\|_{(L^{1,\infty}(w),L^{p})^{\alpha}}:=\sup_{r>0}\ _{r}\left\|f\right\|_{(L^{1,\infty}(w),L^{p})^{\alpha}}.
\end{equation}
We have the following result :% in our weighted spaces. %We denote $(L^{q,\infty}_{w},L^{p})^{\alpha}(\mathbb R^{n})$ the space consisting of measurable functions $f$ such that $\left\|f\right\|_{(L^{q,\infty}_{w},L^{p})^{\alpha}}<\infty$.
%We are now ready to state our main results. 
%Our first result is the boundedness of Calder\'on-Zygmund operators.%$T_{\Omega}$, which is an extension of Theorem 2 of \cite{Wa}.  
%\section{Statment oe main results}
\begin{thm}\label{contT} If $1<q\leq \alpha<p\leq\infty$ and $w\in\mathcal A_{q}$, then the Calder\'on-Zygmund operator $T$ is bounded on $(L^{q}(w),L^{p})^{\alpha}$.

If $q=1$ and $w\in \mathcal A_{1}$, then  we have 
$$\left\| T f\right\|_{(L^{1,\infty}(w),L^{p})^{\alpha}}\lsim \left\|f\right\|_{(L^{q}(w),L^{p})^{\alpha}}.$$ 
\end{thm}
Remark that this result contains Theorem 3.3 in \cite{KS} as a particular case.
% $L^{q,\kappa}(w)$. 
%\begin{thm}[Theorem 3.3 \cite{KS}]\label{rap1}
%If $1<q<\infty$, $0<\kappa<1$ and $w\in\mathcal A_{q}$, then $T$ is bounded on $L^{q,\kappa}_{w}$.
%
%If $q=1$, $0<\kappa<1$ and $w\in\mathcal A_{1}$, then we have %for all $\lambda>0$ and any cube $Q$, we have 
%\begin{equation*}
%\sup_{B:\text{ ball }}\frac{1}{w(B)^{\kappa}}\left\|(Tf)\chi_{B}\right\|^{\ast}_{q_{w},\infty}\lsim\left\|f\right\|_{L^{1,\kappa}_{w}}.
%\end{equation*}  
%\end{thm}

\medskip

In the case $n\geq 2$, we denote by $\mathbb S^{n-1}$ the unit sphere in $\mathbb R^{n}$ equipped with the normalized Lebesgue measure $d\sigma$. For any $\Omega\in L^{\theta}(\mathbb S^{n-1})$ with $1<\theta\leq\infty$, homogeneous of degree zero and such that %with vanishing integral, i.e.,
\begin{equation*}
\int_{\mathbb S^{n-1}}\Omega(x')d\sigma(x')=0,
\end{equation*}
where $x'=x/\left|x\right|$ for any $x\neq 0$, we define the homogeneous singular integral operator $T_{\Omega}$ by
\begin{equation*}
T_{\Omega}f(x)=\mathrm{p.v}\int_{\mathbb R^{n}}\frac{\Omega(y')}{\left|y\right|^{n}}f(x-y)dy,
\end{equation*}
and the Marcinkiewicz integral of higher dimension $\mu_{\Omega}$ by
\begin{equation*}
\mu_{\Omega}(f)(x)=\left(\int^{\infty}_{0}\left|\int_{\left|x-y\right|\leq t}\frac{\Omega(x-y)}{\left|x-y\right|^{n-1}}f(y)dy\right|^{2}\frac{dt}{t^{3}}\right)^{\frac{1}{2}}.
\end{equation*}
%where 
%\begin{equation}
%F_{\Omega,t}(x)=.
%\end{equation}
Duoandikoetxea in \cite{Du} proved that for $\Omega\in L^{\theta}(\mathbb S^{n-1})$ and $1<\theta<\infty$, if $\theta'\leq q<\infty$ and $w\in\mathcal A_{q/\theta'}$ then the operator $T_{\Omega}$ is bounded on $L^{q}(w)$. %One has the following in the weighted Morrey spaces.
%\begin{thm}[Theorem 2 \cite{Wa}]\label{rap2} 
%Assume that $\Omega\in L^{\theta}(\mathbb S^{n-1})$ with $1<\theta<\infty$. Then for every $\theta'\leq q<\infty$, $w\in\mathcal A_{q/\theta'}$ and $0<\kappa<1$, there exists $C>0$ independent of $f$ such that 
%\begin{equation*}
%\left\|T_{\Omega}f\right\|_{L^{q,\kappa}_{w}}\leq C\left\|f\right\|_{L^{q,\kappa}_{w}}.
%\end{equation*}
%\end{thm}
We prove the following extension of Theorem 2 in \cite{Wa}.
\begin{thm}\label{main1}
Let $\Omega\in L^{\theta}(\mathbb S^{n-1})$ with $1<\theta<\infty$. Then for  $\theta'\leq q\leq\alpha< p\leq\infty$ and $w\in\mathcal A_{q/\theta'}$, the operator $T_{\Omega}$ is bounded on $(L^{q}(w),L^{p})^{\alpha}$. %More precisely,
%\begin{equation*}
%\left\|T_{\Omega}(f)\right\|_{q_{w},p,\alpha}\lsim \left\|f\right\|_{q_{w},p,\alpha}.
%\end{equation*}
\end{thm}
As far as Marcinkiewicz operators are concerned, it is proved in \cite{DFP} that if $\Omega\in L^{\theta}(\mathbb S^{n-1})$ and $1<\theta\leq\infty$, then for every  $\theta'< q<\infty$ and $w\in\mathcal A_{q/\theta'}$, there exists $C>0$ such that 
\begin{equation}
\left\|\mu_{\Omega}f\right\|_{L^{q}(w)}\leq C\left\|f\right\|_{L^{q}(w)}.
\end{equation} 
%The corresponding result in weighted Morrey spaces is stated as follows
%\begin{thm}[Theorem 4 \cite{Wa}]\label{rap3}
%Assume that $\Omega\in L^{\theta}(\mathbb S^{n-1})$ with $1<\theta\leq\infty$. Then for every $\theta'< q<\infty$, $w\in\mathcal A_{q/\theta'}$ and $0<\kappa<1$, there exists $C>0$ independent of $f$ such that 
%\begin{equation*}
%\left\|\mu_{\Omega}f\right\|_{L^{q,\kappa}(w)}\leq C\left\|f\right\|_{L^{q,\kappa}_{w}}.
%\end{equation*}
%\end{thm}
\begin{thm}\label{main3}
Let $\Omega\in L^{\theta}(\mathbb S^{n-1})$ with $1<\theta\leq\infty$. Then for  $\theta'< q\leq\alpha< p\leq\infty$ and $w\in\mathcal A_{q/\theta'}$, the operator $\mu_{\Omega}$ is bounded on $(L^{q}(w),L^{p})^{\alpha}$.
\end{thm}
We can find a similar result in Theorem 1.1 of \cite{WT}. %But but we first stated and proved it as one can see in our preprint on arxiv. 
This result for the limit case $p=\infty$, corresponds to Theorem 4 in \cite{Wa}.
\medskip

We also define the Bochner-Riesz operators of order $\delta>0$  %$(n\geq 2)$ are defined initially for Schwartz functions 
in terms of Fourier transforms by 
\begin{equation*}
\left(T^{\delta}_{R}f\right)\widehat{}(\xi)=\left(1-\frac{\left|\xi\right|^{2}}{R^{2}}\right)^{\delta}_{+}\hat{f}(\xi),
\end{equation*}
where $\hat{f}$ denote the Fourier transform of $f$. These operators can be expressed as convolution operators by the formula
\begin{equation}
T^{\delta}_{R}f(x)=(f\ast\phi_{1/R})(x),\label{convbro}
\end{equation}
where $\phi(x)=[(1-\left|\cdot\right|^{2})^{\delta}_{+}]\:\check{}\;(x)$, and $\check{f}$ standing for the inverse Fourier transform of $f$.
The associate maximal operator is defined by
\begin{equation*}
T^{\delta}_{\ast}f(x)=\sup_{R>0}\left|T^{\delta}_{R}f(x)\right|.
\end{equation*} 
Let $n\geq 2$. For $1<q<\infty$ and $w\in\mathcal A_{q}$, Shi and Sun showed in \cite{SS} that $T^{(n-1)/2}_{\ast}$ is bounded on $L^{q}_{w}$. In the limit case $q=1$ we have in \cite{Va} a weak type inequality when $w\in \mathcal A_{1}$, i.e.,% we have in \cite{Va} that there exists $C>0$ such that for all $\lambda>0$, we have
\begin{equation}
\left\|T^{(n-1)/2}_{1}f\right\|_{L^{1,\infty}(w)}\lsim \left\|f\right\|_{L^{1}(w)}.\label{weakbochner}
\end{equation} 
Putting together (\ref{weakbochner}) and the fact that for  a fixed $R>0$
\begin{equation*}
T^{(n-1)/2}_{R}f(x)=(\phi\ast f_{R})_{1/R}(x),
\end{equation*} 
this implies (see \cite{Wa2}) that the weak-type inequality (\ref{weakbochner}) is satisfied for any $R>0$.%, there exists $C>0$ such that for all $\lambda>0$ we have
%\begin{equation}
%w(\left\{x\in\mathbb R^{n}:\left|T^{(n-1)/2}_{R}f(x)\right|>\lambda\right\})\leq\frac{C}{\lambda}\int_{\mathbb R^{n}}\left|f(y)\right|w(y)dy.
%\end{equation}

%\begin{thm}[Theorem 1 \cite{Wa2}]\label{rap4}
%Let $\delta=(n-1)/2$, $1<q<\infty$, $0<\kappa<1$ and $w\in\mathcal A_{q}$. Then there exists $C>0$ such that
%\begin{equation*}
%\left\|T^{\delta}_{\ast}(f)\right\|_{L^{q,\kappa}_{w}}\leq C\left\|f\right\|_{L^{q,\kappa}_{w}}.
%\end{equation*} 
%\end{thm}
\begin{thm}\label{main5}
Let $1\leq q\leq\alpha<p\leq\infty$ and $w\in\mathcal A_{q}$.
\begin{enumerate}
\item If $q>1$ and $\delta=\frac{n-1}{2}$, then $T^{\delta}_{\ast}$ is bounded on $(L^{q}(w),L^{p})^{\alpha}$.\label{as1}
\item If $q=1$ and $\delta=\frac{n-1}{2}$ then for any $R>0$,
\begin{equation*}
\left\|T^{\delta}_{R}f\right\|_{(L^{1,\infty}(w),L^{p})^{^{\alpha}}}\lsim \left\|f\right\|_{(L^{q}(w),L^{p})^{\alpha}}.
\end{equation*}\label{as2}
\end{enumerate}
\end{thm}
%For the operators $T^{\delta}_{R}$, we have the following in the context of Morrey spaces. %weak inequality.
This result contains Theorems 1 and 2 of \cite{Wa2}.
%\begin{thm}[Theorem 2 \cite{Wa2}]\label{rap5}
%Let $\delta=(n-1)/2$, $q=1$, $0<\kappa<1$ and $w\in\mathcal A_{1}$. Then for any given $R>0$,% there exists $C>0$ such that for all balls $B$ 
%\begin{equation*}
%\sup_{B:\text{ ball}}w(B)^{-\kappa}\left\| T^{\delta}_{R}f\chi_{B}\right\|^{\ast}_{1_{w},\infty}\lsim\left\|f\right\|_{L^{1,\kappa}_{w}}.
%\end{equation*}
%\end{thm}

\medskip

We also have results about commutators generated by those operators. %We the same arguments as in the above result, we can obtain the following %generalization of results given in \cite{Wa2} in weighted Morrey space.

We recall that for a linear operator $\mathcal T$ and a locally integrable function $b$,  the commutator operator is defined by
\begin{equation*}
\left[b,\mathcal T\right]f(x)=b(x)\mathcal T f(x)-\mathcal T (bf)(x).
\end{equation*}
In \cite{ST}, it is proved that for the Calder\'on-Zygmund operator $T$ and $b\in BMO$, i.e., the space consisting of locally integrable functions satisfying $\left\|b\right\|_{BMO}<\infty$, where
\begin{equation*}
\left\|b\right\|_{BMO}:=\sup_{B:\text{ ball}}\frac{1}{\left|B\right|}\int_{B}\left|b(x)-b_{B}\right|dx,
\end{equation*} 
with $b_{B}=\frac{1}{\left|B\right|}\int_{B}b(z)dz$, the commutators $\left[b,T\right]$ are bounded in the weighted Lebesgue space $L^{q}(w)$ whenever $1<q<\infty$ and $w\in\mathcal A_{q}$. More precisely, there exists $C>0$ such that
\begin{equation*}
\left\|\left[b,T\right]f\right\|_{L^{q}(w)}\leq C\left\|b\right\|_{BMO}\left\|f\right\|_{L^{q}(w)},
\end{equation*}
for all $f\in L^{q}(w)$. %Noticing that for $w\in\mathcal A_{q}$, $1\leq q<\infty$, \begin{equation*}
%\frac{1}{\left|B\right|}
%\end{equation*}
\begin{thm}\label{contTcom}
Let $b\in BMO$ and $T$ be a Calder\'on-Zygmund operator. If $1<q\leq \alpha<p\leq\infty$, and $w\in\mathcal A_{q}$ then the operator $\left[b,T\right]$ is bounded on $(L^{q}(w),L^{p})^{\alpha}$.
\end{thm}
%For the weighted Morrey spaces we have
%\begin{thm}[Theorem 3.4 \cite{KS}]\label{rap6}
%Let $b\in BMO$ and $T$ be a Calder\'on-Zygmund operator. If $1<q<\infty$, $0<\kappa<1$ and $w\in\mathcal A_{q}$ then $\left[b,T\right]$ is bounded on $L^{q,\kappa}_{w}(\mathbb R^{n})$.
%\end{thm}
For $b\in BMO$, the boundedness of $\left[b,T_{\Omega}\right]$ on $L^{q}(w)$ when $\Omega\in L^{\theta}(\mathbb S^{n-1})$, $1<\theta<\infty$, $\theta'<q<\infty$ and $w\in\mathcal A_{q/\theta'}$ and the one of $\left[b,T^{\delta}_{R}\right]$ on $L^{q}(w)$ for $1<q<\infty$ and $w\in\mathcal A_{q}$ are just consequences of the well-known boundedness criterion for  commutators of linear operators obtained by Alvarez et al in \cite{ABKP}. We deduce from this, the following. 

%In the case of weighted Morrey, the following are proved.
%\begin{thm}[Theorem 3 \cite{Wa}]\label{rap7}
%Assume that $\Omega\in L^{\theta}(\mathbb S^{n-1})$ with $1<\theta<\infty$ and $b\in BMO$. Then for every $\theta'<q<\infty$, $w\in\mathcal A_{q/\theta'}$ and $0<\kappa<1$, the commutator $\left[b,T_{\Omega}\right]$ is bounded on $L^{q,\kappa}_{w}(\mathbb R^{n})$.
%\end{thm} 
%\begin{thm}[Theorem 3 \cite{Wa2}]\label{rap8}
%Let $\delta\geq (n-1)/2$, $1<q<\infty$, $0<\kappa<1$ and $w\in \mathcal A_{q}$. If $b\in BMO$ then for every $R>0$, the operator $\left[b,T^{\delta}_{R}\right]$ is bounded on $L^{q,\kappa}_{w}(\mathbb R^{n})$.
%\end{thm}
%Taking into consideration this result, Komri and Shira\"i proved (see \cite{KS}, Theorem 3.4) the analogue in the weighted Morrey spaces.
%\begin{thm}[Theorem 
\begin{thm}\label{main2}
Let $\Omega\in L^{\theta}(\mathbb S^{n-1})$ with $1<\theta<\infty$ and $b\in BMO$. For every $\theta'<q\leq\alpha<p\leq\infty$ and $w\in\mathcal A_{q/\theta'}$, the commutator $\left[b,T_{\Omega}\right]$ is bounded on $(L^{q}(w),L^{p})^{\alpha}$.
\end{thm}

\begin{thm}\label{main4}
Let $\Omega\in L^{\theta}(\mathbb S^{n-1})$ with $1<\theta<\infty$ and $b\in BMO$. For every $\theta'<q\leq\alpha<p\leq\infty$ and $w\in\mathcal A_{q/\theta'}$, the commutator $\left[b,\mu_{\Omega}\right]$ is bounded on $(L^{q}(w),L^{p})^{\alpha}$.
\end{thm}
The commutators of Marcinkiewicz operators $\mu_{\Omega}$ and a locally integrable function $b$ can be defined by
\begin{equation*}
\left[b,\mu_{\Omega}\right](f)(x)=\left(\int^{\infty}_{0}\left|\int_{\left|x-y\right|\leq t}\frac{\Omega(x-y)}{\left|x-y\right|^{n-1}}[b(x)-b(y)]f(y)dy\right|^{2}\frac{dt}{t^{3}}\right)^{\frac{1}{2}}.
\end{equation*} 
Notice that $\left[b,\mu_{\Omega}\right](f)(x)=\mu_{\Omega}[(b(x)-b)f](x)$.
%where
%\begin{equation}
%F_{\Omega,t}(x)=.
%\end{equation}
For $\Omega\in L^{\theta}(\mathbb S^{n-1})$, $1<\theta\leq\infty$, the boundedness of $\left[b,\mu_{\Omega}\right]$ on $L^{q}(w)$ when $b\in BMO$, $\theta'<q<\infty$ and $w\in\mathcal A_{q/\theta'}$ was established in \cite{DFP}. We have the following result.

% For weighted Morrey, Wang proved
%\begin{thm}[Theorem 5 \cite{Wa}]\label{rap9}
%Assume that $\Omega\in L^{\theta}(\mathbb S^{n-1})$ with $1<\theta\leq \infty$ and $b\in BMO$. Then for every $\theta'<q<\infty$, $w\in\mathcal A_{q/\theta'}$ and $0<\kappa<1$, the operator $\left[b,\mu_{\Omega}\right]$ is bounded on $L^{q,\kappa}_{w}(\mathbb R^{n})$.
%\end{thm}
\begin{thm}\label{combr}
Let $1<q\leq\alpha<p\leq\infty$ and $w\in\mathcal A_{q}$.
 If $\delta\geq\frac{n-1}{2}$ and $b\in BMO$ then the linear commutators $\left[b,T^{\delta}_{R}\right]$ are bounded on $(L^{q}(w),L^{p})^{\alpha}$.
\end{thm}
The case $p=\infty$ in Theorem \ref{contTcom} was proved by Komori and Shirai (see Theorem 3.4 \cite{KS}), while the same case in Theorems \ref{main2}, \ref{main4} and 
\ref{combr} are done by Wang in \cite{Wa} and \cite{Wa2}. 
%\section{Statement of the main results}
% while the next result is giving the boundedness of the integral operator with a rough kernel $\Omega$
%For Marcinkiewicz operators, we have the following extension of Theorem \ref{rap3}.
%
%For Bochner-Riesz operators and their associated maximal functions, we recapitulate in the following the extensions of Theorems \ref{rap4} and \ref{rap5}.% .and \ref{rap8}.
%For commutators, we extend Theorems \ref{rap6}, \ref{rap7}, \ref{rap8} and \ref{rap9}, and obtain the following results.
\section{Proof of the main results}
The following lemma will be the cornerstone in the proofs of our theorems. The results established in \cite{WT} can also be viewed as consequences.
%We will also need the following doubling property of $\mathcal A_{q}$ weights (see Proposition 9.1.5  \cite{Gr}). Let $w\in \mathcal A_{q}$ for some $1<q<\infty$. Then 
%\begin{enumerate}
%for all $\lambda>1$ and all balls $B$, we have 
%\begin{equation}
%w(\lambda B)\lsim \lambda^{nq}w(B).
%\end{equation}
%\item There exists a positive constant $\gamma$ such that we have
%\begin{equation}
%\left(\frac{1}{\left|B\right|}\int_{B}w(t)^{1+\gamma}dt\right)^{\frac{1}{1+\gamma}}\lsim \frac{1}{\left|B\right|}\int_{B}w(t)dt,\text{ for all balls }B,\label{reverseholder}
%\end{equation}
%and for any measurable subset $E$ of a ball $B$, we have
%\begin{equation}
%\frac{w(E)}{w(B)}\lsim \left(\frac{\left|E\right|}{\left|B\right|}\right)^{\frac{\gamma}{1+\gamma}}.\label{contwl}
%\end{equation}
%\end{enumerate} 
%For the proves, we need some lemma which is contained in \cite{Fe3}.
\begin{lem}\label{contgene}
Let $1\leq s\leq q<\infty$, $w\in\mathcal A_{q/s}$ and $\mathcal T:L^{q}_{\mathrm{loc}}(w)\rightarrow L^{q}_{\mathrm{loc}}(w)$ a sub linear operator which satisfies the following property : % Assume that there exist $0\leq\eta<\infty$ such that $\mathcal T$  satisfies 
for all balls $B\subset\mathbb R^{n}$ %and measurable function $g=f\chi_{(2B)^{c}}$
\begin{equation}
\mathcal T(f\chi_{(2B)^{c}})(x)\lsim \sum^{\infty}_{k=1}k\left(\frac{1}{\left|2^{k+1}B\right|}\int_{2^{k+1}B}\left|f(z)\right|^{s}dz\right)^{\frac{1}{s}}\text{ a.e. on }B.\label{hyp1}
\end{equation}
%\left\{\begin{array}{lll}
%\left\|\mathcal T(g)\chi_{B}\right\|_{q_{w}}\lsim \sum^{\infty}_{k=1}(1+\eta k)\left(\frac{w(B)}{w(2^{k+1}B)}\right)^{\frac{1}{q}}\left\|f\chi_{2^{k+1}B}\right\|_{q_{w}}&\text{ if }&q>1\\
%\mathcal T(f\chi_{(2B)^{c}})(x)\lsim \sum^{\infty}_{k=1}(1+\eta k)w(2^{k+1}B)^{-\frac{1}{q}}\left\|f\chi_{2^{k+1}B}\right\|_{q_{w}}%&\text{ if }&q=1\end{array}\right.. 
%for all balls $B$ in $\mathbb R^{n}$ and some .
Then
\begin{enumerate}
\item  if $q>1$ and $\mathcal T$ is bounded on $L^{q}(w)$, then it is also bounded on $(L^{q}(w),L^{p})^{\alpha}$, for $q\leq \alpha<p\leq\infty$,
\item if for all $\lambda>0$
\begin{equation}
w(\left\{x\in\mathbb R^{n}:\left|\mathcal Tf(x)\right|>\lambda\right\})\lsim\frac{1}{\lambda}\int_{\mathbb R^{n}}\left|f(y)\right|w(y)dy,\label{wtbound}
\end{equation}
then for $1\leq\alpha<p\leq\infty$, $\mathcal T$ is bounded from $(L^{1}(w),L^{p})^{\alpha}$ to $(L^{1,\infty}(w),L^{p})^{\alpha}$.
\end{enumerate}
\end{lem}

The proof is partially inspired by \cite{FLY}. The same arguments are used in \cite{DFS} (see also \cite{BFF}), to prove norm inequalities involving Riesz potentials and integral operators satisfying the hypothesis of Theorem 2.1 of \cite{FLY} in the context of $(L^{q},L^{p})^{\alpha}(\mathbb R^{n})$ spaces.
\proof%[Proof of Lemma \ref{contgene}]
Let $1\leq q\leq\alpha<p\leq\infty$ and $f\in (L^{q}(w),L^{p})^{\alpha}$. We fix
%Let $f\in(L^{q}_{w},L^{p})^{\alpha}(\mathbb R^{n})$.  
 $y\in\mathbb R^{n}$ and $r>0$. For almost every $x\in B(y,r)$, we have 
 \begin{eqnarray*}
 \left|\mathcal Tf(x)\right|&\leq&\left|\mathcal T(f\chi_{B(y,2r)})(x)\right|+\left|\mathcal T(f\chi_{B(y,2r)^{c}})(x)\right|\\
 &\lsim&\left|\mathcal T(f\chi_{B(y,2r)})(x)\right|+\sum^{\infty}_{k=1}k\left(\frac{1}{\left|2^{k+1}B\right|}\int_{2^{k+1}B}\left|f(z)\right|^{s}dz\right)^{\frac{1}{s}}
 \end{eqnarray*}
 according to (\ref{hyp1}). 
% \begin{enumerate}
% \item For Write
% \begin{equation*} f(x)=f(x)\chi_{B(y,2r)}(x)+f(x)\chi_{(B(y,2r))^{c}}(x)%\label{decompff}%\sum^{\infty}_{k=1}f(x)\chi_{B(y,2^{k+1}r)\setminus B(y,2^{k}r)}(x)\equiv\sum^{\infty}_{k=1}f_{k}(x)
% \end{equation*}
%\begin{enumerate}
%\item For $q>1$,
% It is sufficient to prove that there exists $C>0$ such that for all $r>0$,
%\begin{equation*}
%\ _{r}\left\|\mathcal Tf\right\|_{q_{w},p,\alpha}\leq C\left\|f\right\|_{q_{w},p,\alpha}.
%\end{equation*}
%  Since $\mathcal T$ is  a sublinear operator, we have
% \begin{equation}
%\left|\mathcal Tf\chi_{B}(x)\right|\leq\left|\mathcal T(f\chi_{2B})\chi_{B}(x)\right|+\left|\mathcal T(f\chi_{(2B)^{c}})\chi_{B}(x)\right|.\label{lem1}
%\end{equation}
\begin{itemize}
\item If $q=s$ then $w\in\mathcal A_{1}$ implies that 
$$\left(\frac{1}{\left|B\right|}\int_{B}\left|f(z)\right|^{q}dz\right)^{\frac{1}{q}}=\left(\frac{1}{w(B)}\int_{B}\frac{w(B)}{\left|B\right|}\left|f(z)\right|^{q}dz\right)^{\frac{1}{q}}\lsim\frac{1}{w(B)^{\frac{1}{q}}}\left\|f\chi_{B}\right\|_{L^{q}(w)}$$
for all balls $B\subset \mathbb R^{n}$.
\item If $s<q$ then H\"older inequality and $\mathcal A_{q/s}$ characterization yield
\begin{eqnarray*}
\left(\frac{1}{\left|B\right|}\int_{B}\left|f(z)\right|^{s}dz\right)^{\frac{1}{s}}&\leq&\left[\frac{1}{\left|B\right|}\left(\int_{B}\left|f(z)\right|^{q}w(z)dz\right)^{\frac{s}{q}}\left(\int_{B}w(z)^{-\frac{s}{q-s}}dz\right)^{\frac{q-s}{q}}\right]^{\frac{1}{s}}\\
&\lsim&\frac{1}{w(B)^{\frac{1}{q}}}\left\|f\chi_{B}\right\|_{L^{q}(w)}.
\end{eqnarray*}
\end{itemize}
It comes that 
\begin{equation}
\left|\mathcal Tf(x)\right|\lsim\left|\mathcal T(f\chi_{B(y,2r)})(x)\right|+\sum^{\infty}_{k=1}\frac{k}{w(B(y,2^{k+1}r))^{\frac{1}{q}}}\left\|f\chi_{B(y,2^{k+1}r)}\right\|_{L^{q}(w)}\label{inpond}
\end{equation}
for almost every $x\in B(y,r)$.
 %\int_{B}\left|f(z)\right|^{q}w(z)dz\right)^{\frac{1}{q}}$$
% \begin{enumerate}
%\item

 First case $q>1$. 

Taking the $L^{q}(w)$-norm on the ball $B(y,r)$ of both sides of (\ref{inpond}) we obtain %The $L^{q}(w)-$boundedness of $\mathcal T$, (\ref{apcondition}) and (\ref{hyp1}) led to %we have %norm of both sides of (\ref{lem1}), we have
 \begin{equation}
 \begin{array}{lll}
\left\|\mathcal Tf\chi_{B(y,r)}\right\|_{L^{q}(w)}&\lsim&\left\|f\chi_{B(y,2r)}\right\|_{L^{q}(w)}\\
&+&\sum^{\infty}_{k=1} k\left\|f\chi_{B(y,2^{k+1}r)}\right\|_{L^{q}(w)}\left(\frac{w(B(y,r))}{w(B(y,2^{k+1}r))}\right)^{\frac{1}{q}}
\end{array}.\label{contr1} 
\end{equation}
according to the boundedness of $\mathcal T$ on $L^{q}(w)$.
Hence, multiplying both sides of (\ref{contr1}) by $w(B)^{\frac{1}{\alpha}-\frac{1}{q}-\frac{1}{p}}$ we obtain %yield
\begin{equation}
\begin{aligned}
%\begin{eqnarray*}
&w(B(y,r))^{\frac{1}{\alpha}-\frac{1}{q}-\frac{1}{p}}\left\|\mathcal Tf\chi_{B(y,r)}\right\|_{L^{q}(w)}\\
&\ \ \ \ \ \ \ \ \ \ \ \ \ \ \ \lsim %w(B(y,2r))^{\frac{1}{\alpha}-\frac{1}{q}-\frac{1}{p}}\left\|f\chi_{B(y,2r)}\right\|_{q_{w}}\\&+&
\sum^{\infty}_{k=0}\frac{k}{2^{kn\frac{s-1}{s}(\frac{1}{\alpha}-\frac{1}{p})}}w(B(y,2^{k+1}r))^{\frac{1}{\alpha}-\frac{1}{q}-\frac{1}{p}}\left\|f\chi_{B(y,2^{k+1}r)}\right\|_{L^{q}(w)},%\left(\frac{w(B)}{w(2^{k+1}B)}\right)^{\frac{1}{\alpha}-\frac{1}{p}}.
\end{aligned}\label{*}
\end{equation}
%\end{eqnarray*}
% while multiplying (\ref{contrw}) by the same quantity with $q=1$ gives 
% \begin{equation*}
%\begin{eqnarray*}
%w(B(y,r))^{\frac{1}{\alpha}-1-\frac{1}{p}}\left\|\mathcal Tf\chi_{B(y,r)}\right\|^{\ast}_{1_{w},\infty}\lsim %w(B(y,2r))^{\frac{1}{\alpha}-\frac{1}{q}-\frac{1}{p}}\left\|f\chi_{B(y,2r)}\right\|_{q_{w}}\\&+&
%\sum^{\infty}_{k=0}\frac{1+\eta k}{2^{\frac{kn}{\tau'}(\frac{1}{\alpha}-\frac{1}{p})}}w(B(y,2^{k+1}r))^{\frac{1}{\alpha}-1-\frac{1}{p}}\left\|f\chi_{B(y,2^{k+1}r)}\right\|_{1_{w}},%\left(\frac{w(B)}{w(2^{k+1}B)}\right)^{\frac{1}{\alpha}-\frac{1}{p}}.
%\end{equation*}
for some $s>1$,  according to Relations (\ref{doubling}) and (\ref{contwl}). Since (\ref{*}) holds for every $y\in\mathbb R^{n}$, this leads to %of the above estimation gives %and having in mind Relations (\ref{doubling}) and (\ref{contwl}), we obtain 
$$\ _{r}\left\|\mathcal Tf\right\|_{(L^{q}(w),L^{p})^{\alpha}}\lsim(1+\sum^{\infty}_{k=1}\frac{ k}{2^{kn\frac{s-1}{s}(\frac{1}{\alpha}-\frac{1}{p})}})\left\|f\right\|_{(L^{q}(w),L^{p})^{\alpha}},\ r>0.$$%&\text{ if }&q>1\\
%\ _{r}\left\|\mathcal Tf\right\|_{(L^{1,\infty}_{w},L^{p})^{\alpha}}\lsim(1+\sum^{\infty}_{k=1}\frac{1+\eta k}{2^{\frac{kn}{\tau'}(\frac{1}{\alpha}-\frac{1}{p})}})\left\|f\right\|_{1_{w},p,\alpha}&\text{ if }&q=1
%\end{array},$$
The expected result follows from taking the supremum over all $r>0$, since  $\sum^{\infty}_{k=1}\frac{k}{2^{kn\frac{s-1}{s}(\frac{1}{\alpha}-\frac{1}{p})}}<\infty$.

Second case $q=1$.

For $\lambda>0$, we have
$$\begin{aligned}
&w(\left\{x\in B(y,r):\left|\mathcal T f(x)\right|>\lambda\right\})\\
&\ \ \ \ \ \ \ \ \ \ \ \ \ \lsim \frac{1}{\lambda}\left(\left\|f\chi_{B(y,2r)}\right\|_{L^{1}(w)}+\sum^{\infty}_{k=1}\frac{kw(B)}{w(B(y,2^{k+1}r))}\left\|f\chi_{B(y,2^{k+1}r)}\right\|_{L^{1}(w)}\right)
\end{aligned}$$
according to (\ref{inpond}). That is, 
\begin{equation*}
\left\|\mathcal T f\chi_{B(y,r)}\right\|_{L^{1,\infty}(w)}\lsim \left\|f\chi_{B(y,2r)}\right\|_{L^{1}(w)}+\sum^{\infty}_{k=1}\frac{kw(B)}{w(B(y,2^{k+1}r))}\left\|f\chi_{B(y,2^{k+1}r)}\right\|_{L^{1}(w)}.
\end{equation*}
Multiplying both sides by $w(B(y,r))^{\frac{1}{\alpha}-1-\frac{1}{p}}$, we conclude as in the case $q>1$. %section 2.% and 
%As for the case $q=1$, the proof is the same except for using Estimation (\ref{wtbound}) instead of the boundedness on $L^{q}_{w}$, Relation (\ref{hyp1}) corresponding to $q=1$ and $w\in\mathcal A_{1}$. This completes the proof.
%\begin{equation}
%\en{equadtion}
%\end{enumerate}
\epf

%$I_{\gamma} $ defined for $0<\gamma<n$ by 
%\begin{equation*}
%I_{\gamma}f(x)=\int_{\mathbb R^{n}}\frac{f(y)}{\left|x-y\right|^{n-\gamma}}dy,
%\end{equation*}
%was bounded from $(L^{q},L^{p})^{^{\alpha}}(\mathbb R^{n})$ to $(L^{q^{\ast}},L^{p})^{^{\alpha^{\ast}}}(\mathbb R^{n})$ where for $s>1$, $\frac{1}{s^{\ast}}=\frac{1}{s}-\frac{\gamma}{n}$.

An immediate application of the above lemma is the following weighted version of Theorem 2.1 in \cite{FLY}.
\begin{prop}\label{corfly}
Let $1< q\leq \alpha<p\leq\infty$. Assume that $\mathcal T$ is a sublinear operator satisfying the property that for any $f\in L^{1}$ with compact support  and $x\notin\mathrm{supp }f $
\begin{equation}
\left|\mathcal Tf(x)\right|\lsim \int_{\mathbb R^{n}}\frac{\left|f(y)\right|}{\left|x-y\right|^{n}}dy.\label{2.1}
\end{equation} 
\begin{enumerate}
\item If for $q>1$ and $w\in\mathcal A_{q}$ the operator $\mathcal T$ is bounded on $L^{q}(w)$ then it is also bounded on $(L^{q}(w),L^{p})^{\alpha}$.
\item If for $w\in \mathcal A_{1}$ we have the weak type estimate 
$$\left\|\mathcal Tf\right\|_{L^{1,\infty}(w)}\lsim \left\|f\right\|_{L^{1}(w)},$$
then we have 
$$\left\|\mathcal Tf\right\|_{(L^{1,\infty}(w),L^{p})^{\alpha}}\lsim \left\|f\right\|_{(L^{1}(w),L^{p})^{\alpha}}.$$
\end{enumerate}
%then $\mathcal T$ is bounded on $(L^{q}_{w},L^{p})^{\alpha}(\mathbb R^{n})$.
\end{prop}
\proof Let $B(y,r)$ be a ball in $\mathbb R^{n}$. For  $x,z\in\mathbb R^{n}$, we have
\begin{equation}
x\in B(y,r)\text{ and }z\notin B(y,2r)\Rightarrow\left|y-z\right|\leq 2\left|z-x\right|\leq 3\left|y-z\right|. \label{estxyz}
\end{equation}
Thus for $x\in B(y,r)$, we have
\begin{eqnarray*}
\left|\mathcal T(f\chi_{(B(y,2r))^{c}})(x)\right|&\lsim&\int_{\mathbb R^{n}}\frac{\left|f\chi_{(B(y,2r))^{c}}(z)\right|}{\left|x-z\right|^{n}}dz
\lsim\sum^{\infty}_{k=1}\int_{2^{k}r\leq\left|y-z\right|<2^{k+1}r}\frac{\left|f(z)\right|}{\left|x-z\right|^{n}}dz\\
&\lsim&\sum^{\infty}_{k=1}\frac{1}{\left|B(y,2^{k+1}r)\right|}\int_{B(y,2^{k+1}r)\setminus B(y,2^{k}r)}\left|f(z)\right|dz.
\end{eqnarray*}
%according to (\ref{estxyz}),% But, H\"older inequality and (\ref{apcondition}) allow to say that
%\begin{equation}
%\int_{B(y,2^{k+1}r)}\left|f(z)\right|dz \lsim\left\|f\chi_{B(y,2^{k+1}r)}\right\|_{q_{w}}\left|B(y,2^{k+1}r)\right|w(B(y,2^{k+1}r))^{-\frac{1}{q}}\label{moyennef}
%\end{equation}
%Therefore, for any ball $B$, we have
%$$\mathcal T(f\chi_{(2B)^{c}})(x)\lsim \sum^{\infty}_{k=1}w(2^{k+1}B)^{-\frac{1}{q}}\left\|f\chi_{2^{k+1}B}\right\|_{q_{w}}$$ 
%for all $x\in B$, 
The conclusion follows from Lemma \ref{contgene}.% We conclude with the the $\mathcal A_{q}$ characterization of $w$.
\epf

%We will often use in the forthcoming proofs, the following implication.
The proof of Theorem \ref{contT} is straightforward from Proposition \ref{corfly}. Notice that this result is valid for non translation invariant CZ operators.

For Theorems \ref{main1}, \ref{main3} and \ref{main5}, we just have to prove that Hypothesis (\ref{hyp1}) of Lemma \ref{contgene} is fulfilled to conclude.% 
%For this purpose, we consider a ball $B=B(y,r)$ centered at some point $y\in\mathbb R^{n}$.  
%We have for $x\in B$% Thus 
%\epf

\proof[Proof of Theorem \ref{main1}]
%Let $f\in(L^{q}_{w},L^{p})^{\alpha}(\mathbb R^{n})$. Since the operator $T_{\Omega}$ is sublinear and is bounded on $L^{q}_{w}$, we just have to prove (\ref{hyp1}).
Let $B=B(y,r)$ be a ball of $\mathbb R^{n}$. %Fix $r>0$ and let $B=B(y,r)$ for some $y\in\mathbb R^{n}$. We put $f=f_{1}+f_{2}$, with $f_{1}=f\chi_{B(y,2r)}$. It is clear that 
%$$\left\|T_{\Omega}f\chi_{B}\right\|_{q_{w}}\leq\left\|T_{\Omega}f_{1}\chi_{B}\right\|_{q_{w}}+\left\|T_{\Omega}f_{2}\chi_{B}\right\|_{q_{w}},$$
%. It comes from the boundedness of $T_{\Omega}$ on the weighted Lebesgue space $L^{q}_{w}$ that 
%$$\left\|T_{\Omega}f_{1}\chi_{B}\right\|_{q_{w}}\lsim \left\|f\chi_{B(y,2r)}\right\|_{q_{w}}.$$
For $x\in B(y,r)$ we have %the term in $f_{2}$ we notice that for all $x\in\mathbb R^{n}$,
\begin{equation*}
\left|T_{\Omega}(f\chi_{(2B)^{c}})(x)\right|\leq\sum^{\infty}_{k=1}\left(\int_{2^{k+1}B\setminus 2^{k}B}\left|\Omega((x-z)')\right|^{\theta}dz\right)^{\frac{1}{\theta}}\left(\int_{2^{k+1}B\setminus 2^{k}B}(\frac{\left|f(z)\right|}{\left|x-z\right|^{n}})^{\theta'}dy\right)^{\frac{1}{\theta'}}
\end{equation*}
by H\"older Inequality. From (\ref{estxyz}), it comes that for $x\in B$ and $z\in 2^{k+1}B\setminus 2^{k}B$, we have $2^{k-1}r\leq \left|x-z\right|\leq 2^{k+2}r$. Thus, for $x\in B(y,r)$ and any positive integer $k$, we have 
\begin{equation}
\left(\int_{2^{k+1}B\setminus 2^{k}B}\left|\Omega((x-z)')\right|^{\theta}dz\right)^{\frac{1}{\theta}}\lsim\left\|\Omega\right\|_{L^{\theta}(\mathbb S^{n-1})}\left|2^{k+1}B\right|^{\frac{1}{\theta}},\label{contomega}
\end{equation}
%Since we also have $\left|x-z\right|\approx\left|y-z\right|$ whenever $x\in B$ and $z\notin 2B$,
and
\begin{equation}
\left(\int_{2^{k+1}B\setminus 2^{k}B}(\frac{\left|f(z)\right|}{\left|x-z\right|^{n}})^{\theta'}dy\right)^{\frac{1}{\theta'}}\lsim\frac{1}{\left|2^{k+1}B\right|}\left(\int_{2^{k+1}B}\left|f(z)\right|^{\theta'}dz\right)^{\frac{1}{\theta'}}.\label{contftheta}
\end{equation}
%according to (\ref{estxyz}). 
Therefore, for any ball $B$ in $\mathbb R^{n}$, we have 
\begin{equation*}
\left|T_{\Omega}(f\chi_{(2B)^{c}})(x)\right|\lsim\sum^{\infty}_{k=1}\left(\frac{1}{\left|2^{k+1}B\right|}\int_{2^{k+1}B}\left|f(z)\right|^{\theta'}dz\right)^{\frac{1}{\theta'}},
\end{equation*}
for all $x\in B$, which ends the proof.% When $\theta'=q$, we have 
%\begin{eqnarray*}
%\left|T_{\Omega}(f\chi_{(2B)^{c}})(x)\right|&\lsim&\sum^{\infty}_{k=1}\left(\frac{1}{w(2^{k+1}B)}\int_{2^{k+1}B}\left|f(z)\right|^{q}w(z)dz\right)^{\frac{1}{q}}\\
%&\lsim&\sum^{\infty}_{k=1}(w(2^{k+1}B))^{-\frac{1}{q}}\left\|f\chi_{2^{k+1}B}\right\|_{q_{w}}
%\end{eqnarray*}
%according to the $\mathcal A_{1}$ characterization of $w$. When $\theta'<q$, let $s=q/\theta'>1$. It follows from H\"older inequality and the $\mathcal A_{s}$ characterization of $w$ that
%\begin{eqnarray*}
%\left|T_{\Omega}(f\chi_{(2B)^{c}})(x)\right|&\lsim&\sum^{\infty}_{k=1}\frac{1}{\left|2^{k+1}B\right|^{\frac{1}{\theta'}}}\left(\int_{2^{k+1}B}\left|f(z)\right|^{q}w(z)dz\right)^{\frac{1}{q}}\left(\int_{2^{k+1}B}w^{-\frac{s'}{s}}(z)dz\right)^{\frac{1}{s'\theta'}}\\
%&\lsim&\sum^{\infty}_{k=1}\left(\frac{1}{w(2^{k+1}B)}\int_{2^{k+1}B}\left|f(z)\right|^{q}w(z)dz\right)^{\frac{1}{q}}.
%\end{eqnarray*}
%Hence for any ball $B$, we have
%\begin{equation*}
%\left|T_{\Omega}(f\chi_{(2B)^{c}})(x)\right|\lsim\sum^{\infty}_{k=1}(w(2^{k+1}B))^{-\frac{1}{q}}\left\|f\chi_{2^{k+1}B}\right\|_{q_{w}}
%\end{equation*}
%for $x\in B$, which end the proof. 
\epf

\proof[Proof of Theorem \ref{main3}]
%Let $f\in (L^{q}_{w},L^{p})^{\alpha}(\mathbb R^{n})$, $1<\theta\leq\infty$ and $\theta'\leq q<\infty$. For $w\in\mathcal A_{q/\theta'}$ it is proved in \cite{DFP} that the sublinear operator $\mu_{\Omega}$ is bounded on $L^{q}_{w}(\mathbb R^{n})$.
 Put $g=f\chi_{(2B)^{c}}$ where $B$ is a ball in $\mathbb R^{n}$. For $x\in B$ and $t>0$ we have 
\begin{equation}
\left\{z:\left|x-z\right|\leq t\right\}\cap(2^{k+1}B\setminus 2^{k}B)\neq\emptyset\Rightarrow t\geq 2^{k-1}r,\label{impl}
\end{equation}
for any positive integer $k$. Therefore, 
\begin{eqnarray*}
\left|\mu_{\Omega}g(x)\right|&=&\left(\int^{\infty}_{0}\left|\int_{(2B)^{c}\cap\left\{z:\left|x-z\right|\leq t\right\}}\frac{\Omega(x-z)}{\left|x-z\right|^{n-1}}f(z)dz\right|^{2}\frac{dt}{t^{3}}\right)^{\frac{1}{2}}\\
&\leq&\sum^{\infty}_{k=1}\left(\int_{2^{k+1}B\setminus 2^{k}B}\frac{\left|\Omega(x-z)\right|}{\left|x-z\right|^{n-1}}f(z)dz\right)\left(\int^{\infty}_{2^{k-1}r}\frac{dt}{t^{3}}\right)^{\frac{1}{2}}\\
&\lsim&\sum^{\infty}_{k=1}\frac{1}{\left|2^{k+1}B\right|^{1/n}}\int_{2^{k+1}B\setminus 2^{k}B}\frac{\left|\Omega(x-z)\right|}{\left|x-z\right|^{n-1}}f(z)dz,
\end{eqnarray*}
for all $x\in B$. We end as in the proof of Theorem \ref{main1}.%Since $\Omega\in L^{\theta}(\mathbb S^{n-1})$, H\"older inequality and the $\mathcal A_{q/\theta'}$ characterization of the weight $w$ led to
%\begin{equation*} \left|\mu_{\Omega}g(x)\right|\lsim\left\|\Omega\right\|_{L^{\theta}}\sum^{\infty}_{k=1}\left(\frac{1}{w(2^{k+1}B)}\int_{2^{k+1}B}\left|f(z)\right|^{q}w(z)dz\right)^{^{\frac{1}{q}}}.
%\end{equation*} 
%for all $x\in B$, and we are done.%We end the proof using Lemma \ref{contgene}.
\epf

\proof[Proof of Theorem \ref{main5}]%Let $f\in (L^{q}_{w},L^{p})^{\alpha}(\mathbb R^{n})$, and $B=B(y,r)$ a ball in $\mathbb R^{n}$. 
As in the above proof, we put $g=f\chi_{(2B)^{c}}$, where $B$ is any ball in $\mathbb R^{n}$.
Since for $R>0$, $T^{\delta}_{R}(g)(x)=\left|g\ast\phi_{1/R}(x)\right|$ with $\phi(x)=\left[(1-\left|\cdot\right|^{2})^{\delta}_{+}\right]\check{}(x)$, we have
\begin{equation*}
\left|(g\ast\phi_{1/R})(x)\right|\leq R^{n}\int_{\mathbb R^{n}}\frac{\left|g(z)\right|}{(R\left|x-z\right|)^{n}}dz=\int_{(2B)^{c}}\frac{\left|f(z)\right|}{\left|x-z\right|^{n}}dz.%\\ &\lsim&\sum^{\infty}_{k=1}\frac{1}{\left|2^{k+1}B\right|}\int_{2^{k+1}B\setminus 2^{k}B}\left|f(z)\right|dz,
\end{equation*}
for $x\in B$, where we use the fact that $\left|\phi(x)\right|\lsim \frac{1}{(1+\left|x\right|)^{\frac{n+1}{2}+\delta}}$ for $\delta\geq\frac{n-1}{2}$.
The Assertions (\ref{as1}) and (\ref{as2}) follow from Proposition \ref{corfly}.
% and equivalence (\ref{estxyz}).
% H\"older inequality and the $\mathcal A_{q}$ characterization of $w$, we have
%\begin{equation}
%\int_{2^{k+1}B}\left|f(z)\right|dz\lsim \frac{\left|2^{k+1}B\right|}{w(2^{k+1}B)^{\frac{1}{q}}}\left\|f\chi_{2^{k+1}B}\right\|_{q_{w}}.\label{moyennef}
%\end{equation}
%Thus we have for all $R>0$,
%\begin{equation}
%\left|T^{\delta}_{R}(g)(x)\right|\lsim \sum^{\infty}_{k=1}w(2^{k+1}B)^{-\frac{1}{q}}\left\|f\chi_{2^{k+1}B}\right\|_{q_{w}}
%\end{equation}
%so that
%\begin{equation}
%\left|T^{\delta}_{\ast}(g)(x)\right|\lsim \sum^{\infty}_{k=1}w(2^{k+1}B)^{-\frac{1}{q}}\left\|f\chi_{2^{k+1}B}\right\|_{q_{w}},
%\end{equation} 
%for all $x\in B$.
\epf

 For the results involving commutators, we need the following properties of  $BMO$ (see \cite{JN}). %(see Theorem 5 of \cite{MW}). Let $w\in\mathcal A_{\infty}$. The norm in the space $BMO(w)$ which is consisted of $b\in L^{1}_{loc}(\mathbb R^{n},w(x)dx)$ such that $\left\|b\right\|_{BMO(w)}<\infty$, where
%\begin{equation*}
%\left\|b\right\|_{BMO(w)}:=\sup_{B:\text{ ball}}\frac{1}{w(B)}\int_{B}\left|b(x)-b_{B,w}\right|w(x)dx
%\end{equation*} 
%and $b_{B,w}=\frac{1}{w(B)}\int_{B}b(x)w(x)dx$, is equivalent to the norm of $BMO(\mathbb R^{n})$.
%We also have the following estimations.
 For $b\in BMO$, $1<q<\infty$ and $w\in \mathcal A_{\infty}$ we have
\begin{equation}
\left\|b\right\|_{BMO}\cong\sup_{B:\text{ ball}}\left(\frac{1}{\left|B\right|}\int_{B}\left|b(x)-b_{B}\right|^{q}dx\right)^{\frac{1}{q}},\label{equivbmo}
\end{equation}
and for all balls $B$
\begin{equation}
\left(\frac{1}{w(B)}\int_{B}\left|b(x)-b_{B}\right|^{q}w(x)dx\right)^{\frac{1}{q}}\lsim \left\|b\right\|_{BMO}.\label{contwbmo}
\end{equation}
Let $b\in BMO$ and $B$ a ball in $\mathbb R^{n}$. For all nonnegative integers $k$, we have
\begin{equation}
\left|b_{2^{k+1}B}-b_{B}\right|\lsim (k+1)\left\|b\right\|_{BMO}. \label{contmoy}
\end{equation}
We also need the following lemma which is just an application of H\"older Inequality, the definition of $\mathcal A_{q}$ weights and Estimation (\ref{contwbmo}). The proof is omitted.
\begin{lem}\label{prodbmof}
Let $1\leq s<q<\infty$. For $b\in BMO$ and $w\in\mathcal A_{q/s}$, we have
\begin{equation*}
\left(\int_{2B\setminus B}\left|b(z)-b_{2B}\right|^{s}\left|f(z)\right|^{s}dz\right)^{\frac{1}{s}}\lsim \left\|b\right\|_{BMO}\left|2B\right|^{\frac{1}{s}}w(2B)^{-\frac{1}{q}}\left\|f\chi_{2B}\right\|_{L^{q}(w)}
\end{equation*}
for all balls $B$ and $f\in L^{q}_{\mathrm{loc}}$.
\end{lem}
Theorems \ref{contTcom}, \ref{main2}, \ref{main4} and \ref{combr} are immediate from the following weighted version of Theorem 2.2 in \cite{FLY}.
\begin{prop}\label{profly}
Let $1\leq\theta<q\leq\alpha<p\leq\infty$ and $w\in\mathcal A_{q/\theta}$. Assume $T$ is a sublinear operator which fulfills conditions (\ref{hyp1}) with $s=\theta$ and admits a commutator with any locally integrable function $b$, satisfying 
\begin{equation}
\left[b,T\right](f)(x)=T[(b(x)-b)f](x).\label{defcomt} %\left|(b(x)-b_{B})T(f)(x)\right|+\left|T((b-b_{B})f)(x)\right|.
\end{equation}
If $\left[b,T\right]$ is bounded on $L^{q}(w)$, then $\left[b,T\right]$ is also bounded on $(L^{q}(w),L^{p})^{\alpha}$.
\end{prop}
\proof
Fix a ball $B=B(y,r)$ in $\mathbb R^{n}$. We have for all $x\in B(y,r)$
\begin{equation*}
\left|\left[b,T\right](f)(x)\right|\lsim \left|\left[b,T\right](f\chi_{2B})(x)\right|+\left|b(x)-b_{B}\right|\left|T(f\chi_{(2B)^{c}})(x)\right|+\left|T[(b_{B}-b)f\chi_{(2B)^{c}}](x)\right|.
\end{equation*}
Thus by the $L^{q}(w)$-boundedness of $\left[b,T\right]$, Relations (\ref{hyp1}) and (\ref{contwbmo}), we have
\begin{eqnarray*}
\left\|\left[b,T\right]f\chi_{B}\right\|_{L^{q}(w)}&\lsim& \left\|f\chi_{2B}\right\|_{L^{q}(w)}
+\left\|b\right\|_{BMO}\sum^{\infty}_{k=1}k\left(\frac{w(B)}{w(2^{k+1}B)}\right)^{\frac{1}{q}}\left\|f\chi_{2^{k+1}B}\right\|_{L^{q}(w)}\\
&+&w(B)^{\frac{1}{q}} \sum^{\infty}_{k=1}k\left(\frac{1}{\left|2^{k+1}B\right|}\int_{2^{k+1}B\setminus2^{k}B}\left|b_{B}-b(z)\right|^{\theta}\left|f(z)\right|^{\theta}dz\right)^{\frac{1}{\theta}}.%\\&=&I+II+III.
\end{eqnarray*}
But then it comes from (\ref{contmoy}) and Lemma \ref{prodbmof} that
\begin{equation*}
\left(\int_{2^{k+1}B}\left|b_{B}-b(z)\right|^{\theta}\left|f(z)\right|^{\theta}dz\right)^{\frac{1}{\theta}}\lsim k\left\|b\right\|_{BMO}\frac{\left|2^{k+1}B\right|^{\frac{1}{\theta}}}{w(2^{k+1}B)^{\frac{1}{q}}}\left\|f\chi_{2^{k+1}B}\right\|_{L^{q}(w)}.%\int_{2^{k+1}B\setminus2^{k}B}\left|f(z)\right|dz+\int_{2^{k+1}B\setminus2^{k}B}\left|b_{2^{k+1}B}-b(z)\right|\left|f(z)\right|dz
\end{equation*}
 Hence, we have
\begin{equation}
\begin{array}{lll}
\left\|\left[b,T\right]f\chi_{B(y,r)}\right\|_{L^{q}(w)}&\lsim& \left\|f\chi_{B(y,2r)}\right\|_{L^{q}(w)}\\
&+&\left\|b\right\|_{BMO}\sum^{\infty}_{k=1} k^{2}\left(\frac{w(B(y,r))}{w(B(y,2^{k+1}r))}\right)^{\frac{1}{q}}\left\|f\chi_{B(y,2^{k+1}r)}\right\|_{L^{q}(w)},
\end{array}\label{comgene}
\end{equation}
for all $y\in\mathbb R^{n}$. Therefore, multiplying both sides of (\ref{comgene}) by $w(B(y,r))^{\frac{1}{\alpha}-\frac{1}{q}-\frac{1}{p}}$, Estimates (\ref{doubling}) and (\ref{contwl}) and the $L^{p}$-norm allow us to obtain
\begin{equation}
\;_{r}\left\|\left[b,T\right]f\right\|_{(L^{q}(w),L^{p})^{\alpha}}\lsim\left\|b\right\|_{BMO}(1+\sum^{\infty}_{k=1}\frac{k^{2}}{2^{nk\frac{s-1}{s}(\frac{1}{\alpha}-\frac{1}{p})}})\left\|f\right\|_{(L^{q}(w),L^{p})^{\alpha}},
\end{equation}
for some constant $s>1$ and all $r>0$. We end the proof by taking the supremum over all $r>0$.
\epf

\end{document}